\documentclass[11pt,leqno]{article}
\usepackage[doc]{optional}
\usepackage{color}
\usepackage{framed}
\usepackage{graphicx}
\definecolor{labelkey}{rgb}{0,0.08,0.45}
\definecolor{refkey}{rgb}{0,0.6,0.0}
\definecolor{Brown}{rgb}{0.45,0.0,0.05}
\definecolor{light}{gray}{.75}
\definecolor{shadecolor}{gray}{.85}%
\usepackage{mathpazo}
\usepackage{amsmath}
\usepackage{amssymb}
\usepackage{theorem}
\newdimen\svparindent
\parindent\svparindent

\newenvironment{svgraybox}%
       {\ClassWarning{Springer-SVMono}{Environment "svgraybox" not available,\MessageBreak
         switching over to "quotation" environment;\MessageBreak
         specify documentclass option "graybox",\MessageBreak
         see SVMono documentation -}%
                \par\addvspace{6pt}
                \list{}{\listparindent12\p@%
                        \leftmargin=12\p@%
                        \itemindent    \listparindent
                        \rightmargin   \leftmargin
                        \parsep        0 \@plus\p@}%
                \expandafter\item\parindent=\svparindent
                \relax\hskip-\listparindent}%
       {\endlist}%

       {\fboxsep=12pt\relax
        \begin{shaded}%
        \list{}{\leftmargin=12pt\rightmargin=2\leftmargin\leftmargin=\z@\topsep=\z@\relax}%
        \list{}{\leftmargin=12pt\rightmargin=2\leftmargin\leftmargin=0\topsep=0\relax}%
        \expandafter\item\parindent=\svparindent
        \hskip-\listparindent}%
       {\endlist\end{shaded}}%

\newcommand{\RXX}{\ensuremath{\,\left[-\infty,+\infty\right]}}
\newcommand{\BB}{\ensuremath{\mathbb B}}
\newcommand{\bZ}[1]{\overleftarrow{\thinspace Z\thinspace}_%
{\negthinspace\negthinspace #1}}
\newcommand{\fr}[1]{\overrightarrow{\thinspace r\thinspace}_%
{\negthinspace\negthinspace #1}}
\newcommand{\fZ}[1]{\overrightarrow{\thinspace Z\thinspace}_%
{\negthinspace\negthinspace #1}}
\newcommand{\br}[1]{\overleftarrow{\thinspace r\thinspace}_%
{\negthinspace\negthinspace #1}}

\newcommand{\thalb}{\ensuremath{\tfrac{1}{2}}}
\newcommand{\menge}[2]{\big\{{#1}~\big |~{#2}\big\}}

\newcommand{\fenv}[1]%
{\ensuremath{\,\overrightarrow{\operatorname{env}}_{#1}}}
\newcommand{\benv}[1]%
{\ensuremath{\,\overleftarrow{\operatorname{env}}_{#1}}}

\newcommand{\scal}[2]{\left\langle{#1},{#2}  \right\rangle}

\newcommand{\RR}{\ensuremath{\mathbb R}}

\newcommand{\RPX}{\ensuremath{\left[0,+\infty\right]}}

\newcommand{\RRX}{\ensuremath{\,\left[-\infty,+\infty\right]}}

\newcommand{\NN}{\ensuremath{\mathbb N}}

\newcommand{\dom}{\ensuremath{\operatorname{dom}}}
\newcommand{\argmin}{\ensuremath{\operatorname{argmin}}}
\newcommand{\argmax}{\ensuremath{\operatorname*{argmax}}}

\newcommand{\intdom}{\ensuremath{\operatorname{int}\operatorname{dom}}\,}

\newcommand{\conv}{\ensuremath{\operatorname{conv}}}

\newcommand{\Id}{\ensuremath{\operatorname{Id}}}

\newcommand{\fproj}[1]{\overrightarrow{P\thinspace}_%
{\negthinspace\negthinspace #1}}
\newcommand{\ffproj}[1]{\overrightarrow{Q\thinspace}_%
{\negthinspace\negthinspace #1}}
\newcommand{\bproj}[1]{\overleftarrow{\thinspace P\thinspace}_%
{\negthinspace\negthinspace #1}}
\newcommand{\bfproj}[1]{\overleftarrow{\thinspace Q\thinspace}_%
{\negthinspace\negthinspace #1}}

\newcommand{\fD}[1]{\overrightarrow{D\thinspace}_%
{\negthinspace\negthinspace #1}}
\newcommand{\ffD}[1]{\overrightarrow{F\thinspace}_%
{\negthinspace\negthinspace #1}}

\newcommand{\bD}[1]{\overleftarrow{\thinspace D\thinspace}_%
{\negthinspace\negthinspace #1}}
\newcommand{\bfD}[1]{\overleftarrow{\thinspace F\thinspace}_%
{\negthinspace\negthinspace #1}}

\newcommand{\pinf}{\ensuremath{+\infty}}

\newcommand{\bc}{\ensuremath{\mathbf{c}}}
\newcommand{\by}{\ensuremath{\mathbf{y}}}

{\begin{list}{}{%
\settowidth{\labelwidth}{\textrm{#1~}}%
\setlength{\leftmargin}{\labelwidth+\labelsep}}}
{\end{list}}
\newtheorem{theorem}{Theorem}[section]
\newtheorem{lemma}[theorem]{Lemma}
\newtheorem{corollary}[theorem]{Corollary}

\newtheorem{definition}[theorem]{Definition}
\theoremstyle{plain}{\theorembodyfont{\rmfamily}
}
\theoremstyle{plain}{\theorembodyfont{\rmfamily}
}
\theoremstyle{plain}{\theorembodyfont{\rmfamily}
}
\theoremstyle{plain}{\theorembodyfont{\rmfamily}
\newtheorem{example}[theorem]{Example}}
\newtheorem{fact}[theorem]{Fact}
\theoremstyle{plain}{\theorembodyfont{\rmfamily}
\newtheorem{remark}[theorem]{Remark}}

\newcommand{\boxedeqn}[1]{%
    \[\fbox{%
        \addtolength{\linewidth}{-2\fboxsep}%
        \addtolength{\linewidth}{-2\fboxrule}%
        \begin{minipage}{\linewidth}%
        \begin{equation}#1\\[+4mm]\end{equation}%
        \end{minipage}%
      }\]%
  }

\newcommand{\boxedregion}[1]{%
\framebox[0.95\linewidth][c]{%
\hspace{2cm}
\addtolength{\linewidth}{-2\fboxsep}%
\addtolength{\linewidth}{-2\fboxrule}%
\begin{minipage}{\linewidth}%
 #1 %
        \end{minipage}%
        \addtolength{\linewidth}{-2\fboxsep}%
        \addtolength{\linewidth}{-2\fboxrule}%
        
      }%
  }

\begin{document}

\title{Chebyshev Sets, Klee Sets, and Chebyshev Centers with
respect to Bregman
Distances:\\ Recent Results and Open Problems}
\author{Heinz H.\ Bauschke, Mason S.\ Macklem, and Xianfu Wang}
\author{
Heinz H.\ Bauschke\thanks{Mathematics, Irving K.\ Barber School,
UBC Okanagan, Kelowna, British Columbia V1V 1V7, Canada. E-mail:
\texttt{heinz.bauschke@ubc.ca}.},
~Mason S.\ Macklem\thanks{Mathematics, Irving K.\ Barber School, UBC Okanagan,
Kelowna, British Columbia V1V 1V7, Canada. E-mail:
\texttt{mason.macklem@ubc.ca}.},
~and Xianfu\
Wang\thanks{Mathematics, Irving K.\ Barber School, UBC Okanagan,
Kelowna, British Columbia V1V 1V7, Canada.
E-mail:  \texttt{shawn.wang@ubc.ca}.}}

%
\maketitle

\def\chaptername{1} 
\def\thechapter{1} 



\abstract{
In Euclidean spaces, the geometric notions
of nearest-points map, farthest-points map,
Chebyshev set, Klee set, and Chebyshev center
are well known and well understood.
Since early works going back to the 1930s, tremendous
theoretical progress has been made, mostly by extending
classical results from Euclidean space to Banach space settings.
In all these results, the distance between points is
induced by some underlying norm.
Recently, these notions have been revisited from a different
viewpoint in which the discrepancy between points is measured by
Bregman distances induced by Legendre functions.
The associated framework covers
the well known Kullback-Leibler divergence and
the Itakura-Saito distance.
In this survey, we review known results and we present new
results on Klee sets and Chebyshev centers with respect to
Bregman distances. Examples are provided and connections to
recent work on Chebyshev functions are made.
We also identify several intriguing open problems.
}

\vspace{0.5cm}

\noindent {\bfseries Keywords:}
Bregman distance,
Chebyshev center, Chebyshev function, Chebyshev point of a function,
Chebyshev set,
convex function,
farthest point,
Fenchel conjugate,
Itakura-Saito distance,
Klee set, Klee function,
Kullback-Leibler divergence,
Legendre function,
nearest point,
projection.

\noindent{\bf AMS 2010 Subject Classification:}
{Primary  41A65;
Secondary 28D05, 41A50, 46N10, 47N10, 49J53, 54E52, 58C06, 90C25.
}

\section{Introduction}

\label{:sec:1}

\subsection*{{Legendre Functions and Bregman Distances}}

Throughout, we assume that
\boxedeqn{
\text{$X=\RR^n$ is
the standard Euclidean space with inner
product $\langle\cdot,\cdot\rangle$,}
}

\noindent with induced norm $\|\cdot\|\colon x\mapsto \sqrt{\scal{x}{x}}$,
and with metric $(x,y)\mapsto \|x-y\|$.
In addition, it is assumed that
\boxedeqn{ 
f\colon X\to\RRX
\text{~is a convex function of Legendre type,}
}
\noindent
also referred as a Legendre function.
We assume the reader is familiar with basic results and standard notation
from Convex Analysis; see, e.g., \cite{Rock70,Rock98,Zali02}.
In particular,
$f^*$ denotes the Fenchel conjugate of $f$, and
$\intdom f$ is the interior of the domain of $f$.
For a subset $C$ of $X$,
$\overline{C}$ stands for the closure of $C$,
$\conv C$ for the convex hull of $C$, and
$\iota_{C}$ for the indicator function of $C$, i.e.,
$\iota_C(x) = 0$, if $x\in C$ and
$\iota_C(x)=+\infty$, if $x\in X\smallsetminus C$.
Now set
\boxedeqn{
U=\intdom f.
}

\begin{example}[Legendre functions]
\label{:ex:Legendre}
\index{Legendre function, Convex function of Legendre type}
The following are Legendre functions
\footnote{Here and elsewhere, inequalities between vectors in
$\RR^n$ are interpreted coordinate-wise.},
each evaluated at a point $x\in X$.
\begin{enumerate}
\item
\label{:ex:Legendre:i}
\emph{halved energy:}
$f(x) = \thalb\|x\|^2 = \thalb\sum_{j}x_j^2$.\\[-1 mm]
\item
\label{:ex:Legendre:ii}
\emph{negative entropy:}
\index{Negative entropy}
$f(x)= \displaystyle \begin{cases}
\sum_{j} \big(x_j\ln(x_j)-x_j\big), &\text{if $x\geq 0$;}\\
\pinf, &\text{otherwise.}
\end{cases}$\\[+2 mm]
\item
\label{:ex:Legendre:iii}
\emph{negative logarithm:}
$f(x)= \displaystyle \begin{cases}
-\sum_{j} \ln(x_j), &\text{if $x>0$;}\\
\pinf, &\text{otherwise.}
\end{cases}$
\end{enumerate}
Note that $U=\RR^n$ in
\ref{:ex:Legendre:i},
whereas $U=\RR^n_{++}$ in \ref{:ex:Legendre:ii}
and \ref{:ex:Legendre:iii}.
\end{example}

Further examples of Legendre functions can be found in,
e.g., \cite{Baus97,commun01,BorVanBook,Rock70}.

\begin{fact}[Rockafellar]
\emph{(See \cite[Theorem~26.5]{Rock70}.)}
\label{:f:Legendre}
The gradient map $\nabla f$ is a continuous
bijection between $\intdom f$ and $\intdom f^*$,
with continuous inverse map $(\nabla f)^{-1}= \nabla f^*$.
Furthermore, $f^*$ is also a convex function of Legendre type.
\end{fact}

Given $x\in U$ and $C\subseteq U$,
it will be convenient to write

\boxedregion{
\begin{align}
\label{:convenient}
x^* &=\nabla f(x),\\
C^* &=\nabla f(C),\\
U^* &=\intdom f^*, \label{:e:a}\\
\nonumber
\end{align}
}

\vspace{0.25cm}
\noindent
and similarly for other vectors and sets in $U$.
Note that we used Fact~\ref{:f:Legendre} for \eqref{:e:a}.

While the Bregman distance defined next
is not a distance in the sense of metric
topology, it does possess some good properties that allow it to
measure the discrepancy between points in $U$.

\begin{definition}[Bregman distance]
\index{Bregman distance}
\emph{(See \cite{Bregman,ButIus,CenZen}.)}
The \emph{Bregman distance} with respect to $f$,
written $D_f$ or simply $D$, is the function
\begin{equation}
D\colon X\times X\to\RPX\colon
(x,y)\mapsto
\begin{cases}
f(x)-f(y) -\scal{\nabla f(y)}{x-y}, &\text{if $y\in U$;}\\
\pinf, &\text{otherwise.}
\end{cases}
\end{equation}
\end{definition}
\noindent

\begin{fact}
\emph{(See \cite[Proposition~3.2.(i) and Theorem~3.7.(iv)\&(v)]{Baus97}.) }
\label{:f:D}
Let $x$ and $y$ be in $U$.
Then the following hold.
\begin{enumerate}
\item
\label{:f:Di}
$D_f(x,y) = f(x)+f^*(y^*)-\scal{y^*}{x} = D_{f^*}(y^*,x^*)$. \\[-1mm]
\item
\label{:f:Dii}
$D_f(x,y)=0$ $\Leftrightarrow$ $x=y$ $\Leftrightarrow$ $x^*=y^*$
$\Leftrightarrow$ $D_{f^*}(x^*,y^*)=0$.
\end{enumerate}
\end{fact}

\begin{example}
\label{:ex:D}
The Bregman distances corresponding to the Legendre functions of
Example~\ref{:ex:Legendre} between two points $x$ and $y$ in $X$ are
as follows.
\begin{enumerate}
\item
\label{:ex:D:i}
$D(x,y) = \tfrac{1}{2}\|x-y\|^2$.\\[+1mm]
\item
\label{:ex:D:ii}
\index{Kullback-Leibler divergence}
$\displaystyle
D(x,y) = \begin{cases}
\sum_j \big(x_j\ln(x_j/y_j)-x_j+y_j\big), &\text{if $x\geq 0$ and $y>0$;}\\
\pinf, &\text{otherwise.}
\end{cases}$ \\[+2mm]
\item
\label{:ex:D:iii}
\index{Itakura-Saito distance}
$\displaystyle
D(x,y) = \begin{cases}
\sum_j \big(\ln(y_j/x_j)+x_j/y_j-1\big), &\text{if $x> 0$ and $y>0$;}\\
\pinf, &\text{otherwise.}
\end{cases}$
\end{enumerate}
These Bregman distances are also known as
\ref{:ex:D:i} the \emph{halved Euclidean distance squared},
\ref{:ex:D:ii} the  \emph{Kullback-Leibler divergence}, and
\ref{:ex:D:iii} the  \emph{Itakura-Saito distance}, respectively.
\end{example}

From now on, we assume that $C$ is a
subset of $X$ such that
\boxedeqn{
\text{$C$ is closed}\quad\text{and}\quad
\varnothing \neq C \subseteq U.
}
\noindent
The \emph{power set} (the set of all subsets) of $C$
is denoted by $2^C$.

We are now in a position
to introduce the various geometric notions.

\subsection*{{Nearest Distance, Nearest Points, and Chebyshev Sets}}

\begin{definition}[Bregman nearest-distance function and
nearest-points map]~\\
\label{:d:DP}
\index{Nearest-distance function, Nearest-points map}
The \emph{left Bregman nearest-distance function} with respect to $C$ is
\begin{equation}
\bD{C}\:  \colon X\to \RPX \colon y\mapsto
\inf_{x\in C}D(x,y),
\end{equation}
and the \emph{left Bregman nearest-points map} with respect to $C$ is
\begin{equation}
\bproj{C}\colon X \to 2^C \colon y\mapsto
\menge{x \in C}{D(x,y) = \bD{C}(y) < \pinf}.
\end{equation}
The \emph{right Bregman nearest-distance} and the
\emph{right Bregman nearest-point map}
with respect to $C$
are
\begin{equation}
\fD{C}\:  \colon X\to \RPX \colon
x\mapsto \inf_{y\in C}D(x,y)
\end{equation}
and
\begin{equation}
\fproj{C}\colon X \to 2^C \colon x\mapsto
\menge{y \in C}{D(x,y) = \fD{C}(x)<\pinf},
\end{equation}
respectively.
If we need to emphasize the underlying Legendre function $f$,
then we write $\bD{f,C}$, $\bproj{f,C}$,
$\fD{f,C}$, and $\fproj{f,C}$.
\end{definition}

\begin{definition}[Chebyshev sets]
\index{Chebyshev set}
The set $C$ is a \emph{left Chebyshev set} with respect to
the Bregman distance,
or simply \emph{$\bD{}$-Chebyshev}, if for every $y\in U$,
$\bproj{C}(y)$ is a singleton.
Similarly, the set $C$ is a \emph{right Chebyshev set}
with respect to the Bregman distance, or
simply \emph{$\fD{}$-Chebyshev}, if for every $x\in U$,
$\fproj{C}(x)$ is a singleton.
\end{definition}

\begin{remark}[Classical Bunt-Motzkin result]
\index{Bunt-Motzkin Theorem}
Assume that $f$ is the halved energy as in
Example~\ref{:ex:Legendre}\ref{:ex:Legendre:i}.
Since the halved Euclidean distance squared (see
Example~\ref{:ex:D}\ref{:ex:D:i}) is symmetric,
the left and right (Bregman) nearest distances coincide, as do
the corresponding nearest-point maps.
Furthermore, the set $C$ is Chebyshev if and only if for every $z\in X$,
the metric\footnote{The metric projection is the nearest-points map
with respect to
the Euclidean distance.} projection $P_C(z)$ is a singleton.
It is well known that if $C$ is convex, then $C$ is Chebyshev.
In the mid-1930s,
Bunt \cite{Bunt} and Motzkin \cite{Motzkin} showed independently that
the following converse holds:
\begin{equation}
\label{:e:BM}
\text{$C$ is Chebyshev $~\Longrightarrow~$ $C$ is convex.}
\end{equation}
For other works in this direction, see, e.g.,
\cite{Edgar1,BerWes,Jon,lewis,Deutsch,urruty1,urruty3,urruty4,Singer70,Singer74,Vlasov73}.
It is still unknown whether or not \eqref{:e:BM} holds in
general Hilbert spaces.
We review corresponding results for the present
Bregman setting in Section~\ref{:sec:ChebSets} below.
\end{remark}

\subsection*{Farthest Distance, Farthest Points, and Klee Sets}

\begin{definition}[Bregman farthest-distance function and
farthest-points map]~\\
\label{:d:FQ}
\index{Farthest-distance function, Farthest-points map}
The \emph{left Bregman farthest-distance function} with respect to $C$ is
\begin{equation}
\bfD{C} \colon X\to \RPX \colon  y\mapsto
\sup_{x\in C}D(x,y),
\end{equation}
and the \emph{left Bregman farthest-points map} with respect to $C$ is
\begin{equation}
\bfproj{C} \colon X \to 2^C \colon
y\mapsto\menge{x\in C}{D(x,y) = \bfD{C}(y)<\pinf}.
\end{equation}
Similarly, the \emph{right Bregman farthest-distance function}
with respect to $C$ is
\begin{equation}
\ffD{C} \colon X \to \RPX \colon  x\mapsto
\sup_{y\in C}D(x,y),
\end{equation}
and the \emph{right Bregman farthest-points map} with respect to $C$ is
\begin{equation}
\ffproj{C} \colon X \to 2^C  \colon
x\mapsto
\menge{y\in C}{D(x,y) =
\ffD{C}(x)<\pinf}.
\end{equation}
If we need to emphasize the underlying Legendre function $f$,
then we write $\bfD{f,C}$, $\bfproj{f,C}$,
$\ffD{f,C}$, and $\ffproj{f,C}$.
\end{definition}

\begin{definition}[Klee sets]
\index{Klee set}
The set $C$ is a \emph{left Klee set } with respect to the
Bregman distance, or simply \emph{$\bD{}$-Klee},
if for every $y\in U$,
$\bfproj{C}(y)$ is a singleton.
Similarly, the set $C$ is a \emph{right Klee set} with respect to the right
Bregman distance, or
simply \emph{$\fD{}$-Klee}, if for every $x\in U$,
$\ffproj{C}(x)$ is a singleton.
\end{definition}

\begin{remark}[Classical Klee result]
\index{Klee Theorem}
Assume again that $f$ is the halved energy as in
Example~\ref{:ex:Legendre}\ref{:ex:Legendre:i}.
Then the left and right (Bregman) farthest-distance functions
coincide, as do
the corresponding farthest-points maps.
Furthermore, the set $C$ is Klee if and only if for every $z\in X$,
the metric farthest-points map
$Q_C(z)$ is a singleton.
It is obvious that if $C$ is a singleton, then $C$ is Klee.
In 1961, Klee \cite{Klee61} showed
the following converse:
\begin{equation}
\label{:e:K}
\text{$C$ is Klee $~\Longrightarrow~$ $C$ is a singleton.}
\end{equation}
See, e.g., also
\cite{Edgar1,lewis,Deutsch,urruty2,urruty3,urruty1,MSV,WesSch}.
Once again, it is still unknown whether or not
\eqref{:e:K} remains true in general Hilbert spaces.
The present Bregman-distance setting is reviewed in
Section~\ref{:sec:Klee} below.
\end{remark}

\subsection*{Chebyshev Radius and Chebyshev Center}

\begin{definition}[Chebyshev radius and Chebyshev center]~\\
\label{:d:rZ}
\index{Chebyshev radius, Chebyshev center}
The left \emph{$\bD{}$-Chebyshev radius} of $C$ is
\begin{equation}
\br{C} = \inf_{y\in U} \bfD{C}(y)
\end{equation}
and the
left \emph{$\bD{}$-Chebyshev center} of $C$ is
\begin{equation}
\bZ{C} = \menge{y\in U}{\bfD{C}(y) = \br{C} <
\pinf}.
\end{equation}
Similarly, the right \emph{$\fD{}$-Chebyshev radius}
of $C$ is
\begin{equation}
\fr{C} = \inf_{x\in U} \ffD{C}(x)
\end{equation}
and the
right \emph{$\fD{}$-Chebyshev center} of $C$ is
\begin{equation}
\fZ{C} = \menge{x\in U}{\ffD{C}(x) = \fr{C} <
\pinf}.
\end{equation}
If we need to emphasize the underlying Legendre function $f$,
then we write $\br{f,C}$, $\bZ{f,C}$,
$\fr{f,C}$, and $\fZ{f,C}$.
\end{definition}

\begin{remark}[Classical Garkavi-Klee result]
\index{Garkavi-Klee Theorem}
Again, assume that $f$ is the halved energy as in
Example~\ref{:ex:Legendre}\ref{:ex:Legendre:i}
so that  the left and right (Bregman) farthest-distance functions
coincide, as do
the corresponding farthest-points maps.
Furthermore, assume that $C$ is bounded.
In the 1960s, Garkavi \cite{Garkavi} and Klee \cite{Klee60} proved
that the Chebyshev center is a singleton,
say $\{z\}$, which is characterized by
\begin{equation}
\label{:e:GK}
z\in\conv Q_C(z).
\end{equation}
See also \cite{NielsenNock,NockNielsen} and
Section~\ref{:sec:ChebCent} below.
In passing, we note that Chebyshev centers are also
utilized in Fixed Point Theory; see, e.g.,
\cite[Chapter~4]{GoebKirk}.
\end{remark}

\subsection*{Goal of the Paper}

The aim of this survey is three-fold.
First, we review recent results
concerning Chebyshev sets, Klee sets, and Chebyshev centers
with respect to Bregman distances.
Secondly, we provide some new results and examples
on Klee sets and Chebyshev centers.
Thirdly, we formulate various tantalizing open problems on
these notions as
well as on the related concepts of Chebyshev functions.

\subsection*{Organization of the Paper}

The remainder of the paper is organized as follows.
In Section~\ref{:sec:aux}, we record auxiliary results
which will make the derivation of the main results more structured.
Chebyshev sets and corresponding open problems are discussed in
Section~\ref{:sec:ChebSets}.
In Section~\ref{:sec:Klee}, we review results and open problems
for Klee sets,
and we also present a new result (Theorem~\ref{:t:18})
concerning left Klee sets.
Chebyshev centers are considered in Section~\ref{:sec:ChebCent},
where we also provide a characterization of left Chebyshev centers
(Theorem~\ref{:t:leftCheb}).
Chebyshev centers are illustrated by Examples in
Section~\ref{:sec:ex}.
{Recent related results on variations of Chebyshev sets and Klee sets are
considered in Section~\ref{:sec:Shawn}.}
Along our journey, we pose several questions that we list collectively
in the final Section~\ref{:sec:op}.

\section{Auxiliary Results}

\label{:sec:aux}

For the reader's convenience, we present the following two results
which are implicitly contained in \cite{bwyy1} and \cite{bwyy2}.

\begin{lemma}
\label{:yuri}
Let $x$ and $y$ be in $C$. Then the following hold.
\begin{enumerate}
\item
$\bD{f,C}(y) = \fD{f^*,C^*}(y^*)$ and
$\fD{f,C}(x) = \bD{f^*,C^*}(x^*)$.\\[-1mm]
\item
$\bproj{f,C}\big|_U =
\nabla f^*\circ \fproj{f^*,C^*}\circ \nabla f$
and
$\fproj{f,C}\big|_U =
\nabla f^*\circ \bproj{f^*,C^*}\circ \nabla f$.\\[-1mm]
\item
$\bproj{f^*,C^*}\big|_{U^*}
= \nabla f \circ \fproj{f,C}\circ \nabla f^*$
and
$\fproj{f^*,C^*}\big|_{U^*}
= \nabla f\circ \bproj{f,C}\circ \nabla f^*$.
\end{enumerate}
\end{lemma}
\begin{proof}
This follows from Fact~\ref{:f:Legendre},
Fact~\ref{:f:D}\ref{:f:Di}, and
Definition~\ref{:d:DP}.
(See also \cite[Proposition~7.1]{bwyy1}.)
\end{proof}

\begin{lemma}
\label{:l:FQ}
Let $x$ and $y$ be in $C$. Then the following hold.
\begin{enumerate}
\item
\label{:l:FQi}
$\bfD{f,C}(y) = \ffD{f^*,C^*}(y^*)$ and
$\ffD{f,C}(x) = \bfD{f^*,C^*}(x^*)$.\\[-1mm]
\item
\label{:l:FQii}
$\bfproj{f,C}\big|_U =
\nabla f^*\circ \ffproj{f^*,C^*}\circ \nabla f$
and
$\ffproj{f,C}\big|_U =
\nabla f^*\circ \bfproj{f^*,C^*}\circ \nabla f$.\\[-1mm]
\item
\label{:l:FQiii}
$\bfproj{f^*,C^*}\big|_{U^*}
= \nabla f \circ \ffproj{f,C}\circ \nabla f^*$
and
$\ffproj{f^*,C^*}\big|_{U^*}
= \nabla f\circ \bfproj{f,C}\circ \nabla f^*$.
\end{enumerate}
\end{lemma}
\begin{proof}
This follows from Fact~\ref{:f:Legendre},
Fact~\ref{:f:D}\ref{:f:Di}, and
Definition~\ref{:d:FQ}.
(See also \cite[Proposition~7.1]{bwyy2}.)
\end{proof}

\bigskip

The next observation on the duality of Chebyshev radii and
Chebyshev centers is new.

\begin{lemma}
\label{:l:rZ}
\index{Chebyshev radius, Chebyshev center}
The following hold.
\begin{enumerate}
\item
\label{:l:rZi}
$\br{f,C} = \fr{f^*,C^*}$ and
$\fr{f,C} = \br{f^*,C^*}$.\\[-1mm]
\item
\label{:l:rZii}
$\bZ{f,C} = \nabla f^*\big( \fZ{f^*,C^*}\big)$ and
$\fZ{f,C} = \nabla f^*\big( \bZ{f^*,C^*}\big)$.\\[-1mm]
\item
\label{:l:rZiii}
$\bZ{f^*,C^*} = \nabla f\big( \fZ{f,C}\big)$ and
$\fZ{f^*,C^*} = \nabla f\big( \bZ{f,C}\big)$.\\[-1mm]
\item
\label{:l:rZiv}
$\bZ{f,C}$ is a singleton $\Leftrightarrow$
$\fZ{f^*,C^*}$ is a singleton.\\[-1mm]
\item
\label{:l:rZv}
$\fZ{f,C}$ is a singleton $\Leftrightarrow$
$\bZ{f^*,C^*}$ is a singleton.
\end{enumerate}
\end{lemma}
\begin{proof}
\ref{:l:rZi}:
Using Definition~\ref{:d:rZ}
and Lemma~\ref{:l:FQ}\ref{:l:FQi}, we see that
\begin{equation}
\br{f,C} = \inf_{y\in U} \bfD{C}(y) =
\inf_{y^*\in U^*} \ffD{C^*}(y^*) =
\fr{f^*,C^*}
\end{equation}
and that
\begin{equation}
\fr{f,C} = \inf_{y\in U} \ffD{C}(y) =
\inf_{y^*\in U^*} \bfD{C^*}(y^*) =
\br{f^*,C^*}.
\end{equation}
\ref{:l:rZii}\&\ref{:l:rZiii}:
Let $z\in U$.
Using \ref{:l:rZi} and Lemma~\ref{:l:FQ}\ref{:l:FQi},
we see that
\begin{equation}
z\in\bZ{f,C} \Leftrightarrow
\bfD{f,C}(z)=\br{f,C}
\Leftrightarrow
\ffD{f^*,C^*}(z^*)=\fr{f^*,C^*}
\Leftrightarrow
z^*\in\fZ{f^*,C^*}.
\end{equation}
This verifies $\bZ{f,C} = \nabla f^*\big( \fZ{f^*,C^*}\big)$
and $\fZ{f^*,C^*} = \nabla f\big( \bZ{f,C}\big)$.
The remaining identities follow similarly.

\ref{:l:rZiv}\&\ref{:l:rZv}:
Clear from \ref{:l:rZii}\&\ref{:l:rZiii}
and Fact~\ref{:f:Legendre}.
\end{proof}

The following two results play a key role for studying the single-valuedness of
$\fproj{f,C}$ via $\bproj{f^*,C^*}$ and $\ffproj{f,C}$
via $\bfproj{f^*,C^*}$ by duality.

\begin{lemma}\label{:onecategory}
Let $ V$ and $W$ be nonempty open subsets of $X$,
and let $T:V\to W$ be a homeomorphism, i.e., $T$ is a bijection and
both $T$ and $T^{-1}$ are continuous.
Furthermore, let $G$ be a residual\,\footnote{also known as
``second category''} subset of $V$.
Then $T(G)$ is a residual subset of $W$.
\end{lemma}
\begin{proof}
As $G$ is residual,
there exists a sequence of dense open subsets $(O_{k})_{k\in\NN}$
of $V$ such that $G\supseteq \bigcap_{k\in\NN}O_{k}$.
Then $T(G)\supseteq T(\bigcap_{k\in\NN}O_{k})=
\bigcap_{k\in\NN}T(O_{k})$.
Since $T:V\to W$ is a homeomorphism and each $O_{k}$
is dense in $V$, we see that each $T(O_{k})$ is open and
dense in $W$. Therefore,
$\bigcap_{k\in\NN}T(O_{k})$ is a dense $G_{\delta}$ subset in $W$.
\end{proof}

\begin{lemma}
\label{:twomeasure}
Let $V$ be a nonempty open subset of $X$, and let
$T: V\to\RR^n$ be locally Lipschitz.
Furthermore, let $S$ be a subset of $V$ that has Lebesgue measure zero.
Then $T(S)$ has Lebesgue measure zero as well.
\end{lemma}

\begin{proof}
Denote the closed unit ball in $X$ by $\BB$.
For every $y\in V$, let $r(y)>0$ be such that
$T$ is Lipschitz continuous with constant $c(y)$ on the
open ball $O(y)$ centered at $y$ of radius $r(y)$.
In this proof we denote the Lebesgue measure by $\lambda$.
Let $K$ be a compact subset of $X$. To show that $T(S)$ has Lebesgue
measure zero, it suffices to show that $\lambda(T(K\cap S))=0$ because
\begin{equation}
\lambda\big(T(S)\big)=
\lambda\bigg(T\Big(\bigcup_{k\in\NN}S \cap k\BB\Big)
\bigg)\leq
\sum_{k\in\NN}\lambda\big( T(S \cap k\BB)\big).
\end{equation}
The Heine-Borel theorem provides a finite subset
$\{y_{1},\ldots, y_{m}\}$ of $V$
such that
\begin{equation}
\label{:e:100315}
K\subseteq\bigcup_{j=1}^{m}O(y_{j}).
\end{equation}
We now proceed using a technique implicit in
the proof of \cite[Corollary~1]{guzman}.
Set $c=\max\{c_{1},c_{2},\ldots, c_{m}\}$.
Given $\varepsilon >0$, there exists an
open subset $G$ of $X$ such that $G\supseteq K\cap S$
and $\lambda(G)<\varepsilon$.
For each $y\in K\cap S$, let $Q(y)$ be an open cubic interval
centered at $y$ of semi-edge length $s(y)>0$ such that
\begin{equation}\label{:timechange1}
(\exists\,j\in\{1,\ldots,m\})\quad
Q(y)\subseteq G\cap O(y_{j}).
\end{equation}
Then for each $x\in Q(y)$, we have
\begin{equation}
\|Tx-Ty\|\leq c\|x-y\|\leq c\sqrt{n}s(y).
\end{equation}
Hence the image of $Q(y)$ by $T$, $T(Q(y))$,
is contained in a cubic interval --- which we denote by $Q^*(Ty)$ ---
of center $Ty$ and with semi-edge length $c\sqrt{n}s(y)$.
Applying the Besicovitch Covering Theorem,
we see that there exists a sequence $(Q_{k})_{k\in \NN}$
chosen among the open covering $(Q(y))_{y\in K\cap S}$
such that
\begin{equation}\label{:timechange2}
K\cap S\subseteq \bigcup_{k\in\NN}Q_{k}
\quad \text{ and }\quad \sum_{k\in\NN}\chi_{Q_{k}}\leq \theta,
\end{equation}
where $\chi_{Q_k}$ stands for the characteristic function of $Q_k$ and
where the constant $\theta$ only depends on the dimension of $X$.
Thus,
\begin{equation}
T(K\cap S)\subseteq T\Big(\bigcup_{k\in\NN}Q_{k}\Big)
=\bigcup_{k\in\NN}T(Q_{k})\subseteq
\bigcup_{k\in\NN}Q_{k}^*.
\end{equation}
Now set $d = ( c\sqrt{n})^n$
so that $\lambda(Q_{k}^*)\leq d \lambda(Q_{k})$.
Then, using \eqref{:timechange1} and \eqref{:timechange2},
we see that
\begin{align}
\lambda\Big(\cup_{k\in\NN}Q_{k}^*\Big)
&\leq \sum_{k\in\NN}\lambda(Q_{k}^*)\leq
d\sum_{k\in\NN}\lambda(Q_{k})
=d\sum_{k\in\NN}\int \chi_{Q_{k}}
=d\int \sum_{k\in\NN}\chi_{Q_{k}}\leq
d\theta\lambda(G)\notag\\
&\leq d\theta\varepsilon.
\end{align}
Since $\varepsilon$ was chosen arbitrarily,
we conclude that $\lambda(T(K\cap S))=0$.

Alternatively, one may argue as follows starting from
\eqref{:e:100315}.
We have
$K\cap S\subseteq \big(\bigcup_{j=1}^{m}O(y_{j})\big)\cap S
=\bigcup_{j=1}^{m}O(y_{j})\cap S$
so that
\begin{equation}
\label{:royden}
T(K\cap S)\subseteq \bigcup_{j=1}^{m}T(O(y_{j})\cap S).
\end{equation}
Since $T$ is Lipschitz on each $O(y_{j})$ with constant
$c(y_{j})$ and since $\lambda(O(y_{j})\cap S)=0$,
we apply \cite[Proposition~262D, page~286]{fremlin} and conclude
that $\lambda(T(O(y_{j})\cap S))=0$.
Therefore, $\lambda(T(K\cap S))=0$ by \eqref{:royden}.
\end{proof}

\section{Chebyshev Sets}

\label{:sec:ChebSets}

We start by reviewing the strongest known results
concerning left and right
Chebyshev sets with respect to Bregman distances.

\index{Chebyshev set}

\begin{fact}[$\bD{}$-Chebyshev sets]
\label{:f:leftCheb}
\emph{(See \cite[Theorem~4.7]{bwyy1}.)}
Suppose that $f$ is supercoercive\footnote{
By \cite[Proposition~2.16]{Baus97}, $f$ is supercoercive $:\Leftrightarrow$
$\displaystyle \lim_{\|x\|\to\pinf} \frac{f(x)}{\|x\|} =
\pinf$ $\Leftrightarrow$ $\dom f^* = X$.}
and that $C$ is $\bD{}$-Chebyshev.
Then $C$ is convex.
\end{fact}

\begin{fact}[$\fD{}$-Chebyshev sets]
\label{:f:rightCheb}
\emph{(See \cite[Theorem~7.3]{bwyy1}.)}
Suppose that $\dom f = X$,
that $\overline{C^*}\subseteq U^*$, and
that $C$ is $\fD{}$-Chebyshev.
Then $C^*$ is convex.
\end{fact}

It is not known whether or not Fact~\ref{:f:leftCheb}
and Fact~\ref{:f:rightCheb} are the best possible results.
For instance, is the assumption on supercoercivity in
Fact~\ref{:f:leftCheb} really necessarily?
Similarly, do we really require full domain of $f$ in
Fact~\ref{:f:rightCheb}?

\begin{example}
\label{:ex:curious}
{(See \cite[Example~7.5]{bwyy1}.)}
Suppose that $X=\RR^2$,
that $f$ is the negative entropy
\index{Negative entropy, Kullback-Leibler divergence}
(see Example~\ref{:ex:Legendre}\ref{:ex:Legendre:ii}),
and that
\begin{equation}
C = \menge{(e^\lambda,e^{2\lambda})}{\lambda\in[0,1]}.
\end{equation}
Then $f$ is supercoercive and
$C$ is a \emph{nonconvex}  $\fD{}$-Chebyshev set.
\end{example}

Example~\ref{:ex:curious} is somewhat curious ---
not only does it illustrate that the right-Chebyshev-set counterpart
of Fact~\ref{:f:leftCheb} fails but it also shows
that the conclusion of Fact~\ref{:f:rightCheb} may hold even though
$f$ is not assumed to have full domain.

\begin{fact}
\label{:f:Noll}
\emph{(See \cite[Lemma~3.5]{BN}.)}
Suppose that $f$ is the negative entropy
(see Example~\ref{:ex:Legendre}\ref{:ex:Legendre:ii})
and that $C$ is convex.
Then  $C$ is $\fD{}$-Chebyshev.
\end{fact}

Fact~\ref{:f:Noll} raises two intriguing questions.
Apart from the case of quadratic functions,
are there instances of $f$ where $f$ has full domain and
where every closed convex subset of $U$ is
$\fD{}$-Chebyshev? Because of Fact~\ref{:f:rightCheb},
an affirmative answer to this question
would imply that
$\nabla f$ is
a (quite surprising)
\emph{nonaffine yet convexity-preserving} transformation.
Combining Example~\ref{:ex:curious} and Fact~\ref{:f:Noll}, we
deduce that --- when working with the negative entropy ---
if $C$ is convex, then $C$ is $\fD{}$-Chebyshev
but \emph{not} vice versa.
Is it possible to describe the $\fD{}$-Chebyshev sets in this
setting?

We also note that $C$ is ``nearly $\bD{}$-Chebyshev''
in the following sense.

\begin{fact}
\emph{(See \cite[Corollary~5.6]{bwyy1}.)}
\label{:f:smallbP}
Suppose that $f$ is supercoercive,
that $f$ is twice continuously differentiable, and
that for every $y\in U$, $\nabla^2 f(y)$ is positive definite.
Then $\bproj{C}$
is almost everywhere and generically\footnote{That is,
the set $S$
of points $y\in U$ where $\bproj{C}(y)$ is \emph{not}
a singleton is very small both in measure theory
($S$ has measure $0$) and in category theory
($S$ is meager/first category).}
single-valued on $U$.
\end{fact}

It would be interesting to see whether or not supercoercivity is
essential in Fact~\ref{:f:smallbP}.
By duality, we obtain the following result on the
single-valuedness of  $\fproj{f,C}$.

\begin{corollary}
\label{:c:sushi}
Suppose that $f$ has full domain,
that $f^*$ is twice continuously differentiable,
and that $\nabla^2 f^*(y)$ is positive definite for every $y\in U^*$.
Then $\fproj{f,C}$ is almost everywhere and generically
single-valued on $U$.
\end{corollary}
\begin{proof}
By Lemma~\ref{:yuri}(ii),
$\fproj{f,C}\big|_U =
\nabla f^*\circ \bproj{f^*,C^*}\circ \nabla f$.
Fact~\ref{:f:smallbP} states that
$\bproj{f^*,C^*}$ is almost everywhere and generically
single-valued on $U^*$.
Since $f^*$ is twice continuously differentiable, it follows from the Mean
Value Theorem that $\nabla f^*$ is locally Lipschitz.
Since $(\nabla f)^{-1}=\nabla f^*$ is a locally Lipschitz homeomorphism
from $U^*$ to $U$, the conclusion follows from
Lemma~\ref{:onecategory} and Lemma~\ref{:twomeasure}.
\end{proof}

\section{Klee Sets}

\label{:sec:Klee}

\index{Klee set}

Previously known were the following two results.

\begin{fact}[$\bD{}$-Klee sets]
\label{:f:leftKlee}
\emph{(See \cite[Theorem~4.4]{bwyy2}.)}
Suppose that $f$ is supercoercive, that $C$ is bounded, and that
$C$ is $\bD{}$-Klee.
Then $C$ is a singleton.
\end{fact}

\begin{fact}[$\fD{}$-Klee sets]
\label{:f:rightKlee}
\emph{(See \cite[Theorem~3.2]{BMSW}.)}
Suppose that $C$ is bounded and that
$C$ is $\fD{}$-Klee.  Then $C$ is a singleton.
\end{fact}

Fact~\ref{:f:leftKlee} immediately raises the question on whether or
not supercoercivity is really an essential hypothesis.
Fortunately, thanks to Fact~\ref{:f:rightKlee}, which was recently
proved for general Legendre functions without any further assumptions, we
are now able to present a new result which
removes the supercoercivity assumption in Fact~\ref{:f:leftKlee}.

\begin{theorem}[$\bD{}$-Klee sets revisited]
\label{:t:18}
Suppose that $C$ is bounded and that $C$ is $\bD{}$-Klee.
Then $C$ is a singleton.
\end{theorem}
\begin{proof}
On the one hand, since $C$ is compact,
Fact~\ref{:f:Legendre} implies that $C^*$ is compact.
On the other hand, by Lemma~\ref{:l:FQ}\ref{:l:FQiii},
the set $C^*$ is $\fD{f^*}$-Klee.
Altogether, we deduce from Fact~\ref{:f:rightKlee} (applied to $f^*$
and $C^*$) that $C^*$ is a singleton.
Therefore, $C$ is a singleton by Fact~\ref{:f:Legendre}.
\end{proof}

Similarly to the setting of Chebyshev sets,
the set $C$ is ``nearly $\bD{}$-Klee''
in the following sense.

\begin{fact}
\emph{(See \cite[Corollary~5.2.(ii)]{bwyy1}.)}
\label{:f:smallbQ}
Suppose that $f$ is supercoercive,
that $f$ is twice continuously differentiable,
that for every $y\in U$, $\nabla^2 f(y)$ is positive definite,
and that $C$ is bounded.
Then $\bfproj{C}$
is almost everywhere and generically single-valued on $U$.
\end{fact}

Again, it would be interesting to see whether or not supercoercivity is
essential in Fact~\ref{:f:smallbQ}.
Similarly to the proof of Corollary~\ref{:c:sushi},
we obtain the following result on
the single-valuedness of $\ffproj{f,C}$.

\begin{corollary}
Suppose that $f$ has full domain,
that $f^*$ is twice continuously differentiable,
that $\nabla^2 f^*(y)$ is positive definite for every $y\in U^*$,
and that $C$ is bounded.
Then
$\ffproj{f,C}$ is almost everywhere and generically
 single-valued on $U$.
\end{corollary}

\section{Chebyshev Centers: Uniqueness and Characterization}

\label{:sec:ChebCent}

\index{Chebyshev center}

\begin{fact}[$\fD{}$-Chebyshev centers]
\label{:f:rightChebcent}
\emph{(See \cite[Theorem~4.4]{BMSW}.)}
Suppose that $C$ is bounded.
Then the right Chebyshev center with respect to $C$ is a singleton, say
$\fZ{C} = \{x\}$, and $x$ is characterized by
\begin{equation}
x\in\nabla f^*\big(\conv \nabla f(\ffproj{C}(x)) \big).
\end{equation}
\end{fact}

We now present a corresponding new result on the left Chebyshev center.

\begin{theorem}[$\bD{}$-Chebyshev centers]
\label{:t:leftCheb}
Suppose that $C$ is bounded.
Then the left Chebyshev center with respect to $C$ is a singleton, say
$\bZ{C} = \{y\}$, and $y$ is characterized by
\begin{equation}
\label{:e:kati3}
y\in\conv \bfproj{C}(y).
\end{equation}
\end{theorem}
\begin{proof}
By Lemma~\ref{:l:rZ}\ref{:l:rZii},
\begin{equation}
\label{:e:kati2}
\bZ{f,C} = \nabla f^*\big( \fZ{f^*,C^*} \big).
\end{equation}
Now $C^*$ is a bounded subset of $U^*$ because
of the compactness of $C$ and Fact~\ref{:f:Legendre}.
Applying Fact~\ref{:f:rightChebcent} to $f^*$ and $C^*$, we obtain
that $\fZ{f^*,C^*} = \{y^*\}$ for some $y^*\in U^*$ and that
$y^*$ is characterized by
\begin{equation}
\label{:e:kati1}
y^*\in\nabla f\big(\conv \nabla f^*(\ffproj{f^*,C^*}(y^*))
\big).
\end{equation}
By \eqref{:e:kati2},
$\bZ{f,C} = \nabla f^*\big( \fZ{f^*,C^*} \big)
= \{\nabla f^*(y^*)\} = \{y\}$ is a singleton.
Moreover, using Lemma~\ref{:l:FQ}\ref{:l:FQii},
we see that the characterization \eqref{:e:kati1} becomes
\begin{align}
\bZ{f,C} = \{y\}
&\Leftrightarrow y^*\in\nabla f\big(\conv \nabla f^*(\ffproj{f^*,C^*}(y^*))
\big)\notag\\
&\Leftrightarrow \nabla f^*(y^*)\in\conv \nabla
f^*(\ffproj{f^*,C^*}(y^*))\notag\\
&\Leftrightarrow y\in\conv \nabla
f^*(\ffproj{f^*,C^*}(\nabla f(y)))\notag\\
&\Leftrightarrow y\in\conv
\bfproj{f,C}(y),
\end{align}
as claimed.
\end{proof}

\begin{remark}
The proof of Fact~\ref{:f:rightChebcent} does not carry over
directly to the setting of Theorem~\ref{:t:leftCheb}.
Indeed, one key element in that proof was to realize that
the right farthest distance function
\begin{equation}
\ffD{C} = \sup_{y\in C} D(\cdot,y)
\end{equation}
is \emph{convex} (as the supremum of convex functions) and
then to apply the Ioffe-Tihomirov theorem (see, e.g.,
\cite[Theorem~2.4.18]{Zali02}) for the subdifferential of the
supremum of convex function.
In contrast,
$\bfD{C} = \sup_{x\in C} D(x,\cdot)$
is generally \emph{not convex}.
(For more on separate and joint convexity of $D$, see \cite{BB01}.)
\end{remark}

\section{Chebyshev Centers: Two Examples}

\index{Chebyshev center}

\label{:sec:ex}

\subsection*{Diagonal-Symmetric Line Segments in the Strictly
Positive Orthant}

In addition to our standing assumptions from Section~\ref{:sec:1},
we assume in this Subsection that the following hold:

\boxedregion{
\begin{gather}
X = \mathbb{R}^2;\\
\bc_0 = (1,a)\;\;\text{and}\;\;\bc_1=(a,1),\quad
\text{where}\;\; 1<a<\pinf;\\
\bc_\lambda = (1-\lambda)\bc_0 + \lambda \bc_1,\quad
\text{where}\;\; 0<\lambda<1;\\
C = \conv\big\{\bc_0,\bc_1\big\}
=\menge{\bc_\lambda}{0\leq\lambda\leq 1}.\\
\nonumber
\end{gather}
}



\begin{theorem}
\label{:t:segment}
Suppose that $f$ is any of the functions considered in
Example~\ref{:ex:Legendre}.
Then the left Chebyshev center is the midpoint of $C$, i.e.,
$\bZ{C} = \{\bc_{1/2}\}$.
\end{theorem}
\begin{proof}
By Theorem~\ref{:t:leftCheb}, we write
$\bZ{C}=\{\by\}$, where $\by=(y_1,y_2)\in U$.
In view of \eqref{:e:kati3} and
Fact~\ref{:f:D}\ref{:f:Dii}, we obtain that
$\bfproj{C}(\by)$ contains at least two elements.
On the other hand, since $\bfproj{C}(\by)$ consists
of the maximizers of the \emph{convex} function $D(\cdot,\by)$ over
the compact set $C$, \cite[Corollary~32.3.2]{Rock70} implies that
$\bfproj{C}(\by)\subseteq\{\bc_0,\bc_1\}$.
Altogether,
\begin{equation}
\label{:e:kati6}
\bfproj{C}(\by)=\big\{\bc_0,\bc_1\big\}.
\end{equation}
In view of \eqref{:e:kati3},
\begin{equation}
\label{:e:kati4}
\by \in C.
\end{equation}
On the other hand, a symmetry argument identical to the proof of
\cite[Proposition~5.1]{BMSW} and the uniqueness of its Chebyshev center show that
$\by$ must lie
on the diagonal, i.e., that
\begin{equation}
\label{:e:kati5}
y_1 = y_2.
\end{equation}
The result now follows because the only point satisfying both
\eqref{:e:kati4} and \eqref{:e:kati5} is
$\bc_{1/2}$, the midpoint of $C$.
\end{proof}

\begin{remark}
Theorem~\ref{:t:segment} is in stark contrast
with \cite[Section~5]{BMSW}, where
we investigated the \emph{right} Chebyshev center in this setting.
Indeed, there we found that the right Chebyshev center does depend
on the underlying Legendre function used
(see \cite[Examples~5.2, 5.3, and 5.5]{BMSW}).
Furthermore, for each Legendre function $f$
considered in Example~\ref{:ex:Legendre}, we obtain the following
formula.
\begin{equation}
\big(\forall \by=(y_1,y_2)\in U\big)\quad
\bfproj{f,C}(\by) =
\begin{cases}
\{\bc_0\}, &\text{if $y_2<y_1$;}\\
\{\bc_1\}, &\text{if $y_2>y_1$;}\\
\{\bc_0,\bc_1\}, &\text{if $y_1=y_2$.}
\end{cases}
\end{equation}
Indeed, since for every $\by\in U$,
the function $D(\cdot,\by)$ is convex;
the points where the supremum is achieved is a subset of
the extreme points of $C$, i.e., of $\{\bc_{0}, \bc_{1}\}$.
Therefore, it suffices to compare $D(\bc_{0}, \by)$
and $D(\bc_{1},\by)$.
\end{remark}

\subsection*{Intervals of Real Numbers}

\begin{theorem}
\label{:t:kuehl}
Suppose that $X=\RR$ and that $C=[a,b]\subset U$, where
$a\neq b$.
Denote the right and left Chebyshev centers by $x$ and $y$, respectively.
Then\,\footnote{Recall the convenient notation introduced
on page~\pageref{:convenient}!}
\begin{equation}
x = \frac{f^*(b^*)-f^*(a^*)}{b^*-a^*}
\quad\text{and}\quad
y^* = \frac{f(b)-f(a)}{b-a}.
\end{equation}
\end{theorem}
\begin{proof}
Analogously to the derivation of
\eqref{:e:kati6}, it must hold that
\begin{equation}
\bfproj{C}(y) = \{a,b\}.
\end{equation}
This implies that $y$ satisfies
$D(a,y) = D(b,y)$. In turn, using Fact~\ref{:f:D}\ref{:f:Di},
this last equation is equivalent to
$D_{f^*}(y^*,a^*) = D_{f^*}(y^*,b^*)$
$\Leftrightarrow$
$f^*(y^*)+f(a)-y^*a = f^*(y^*) + f(b)-y^*b$
$\Leftrightarrow$
$f(b)-f(a) = y^*(b-a)$
$\Leftrightarrow$
$y^* = (f(b)-f(a))/(b-a)$, as claimed.
Hence
\begin{equation}
y = \nabla f^*\Big( \frac{f(b)-f(a)}{b-a}\Big).
\end{equation}
Combining this formula (applied to $f^*$ and $C^*=[a^*,b^*]$) with
Lemma~\ref{:l:rZ}\ref{:l:rZii}, we obtain that the right
Chebyshev center is given by
\begin{equation}
x = \nabla f^*\left( \nabla
f^{**}\Big(\frac{f^*(b^*)-f^*(a^*)}{b^*-a^*}\Big)\right)
= \frac{f^*(b^*)-f^*(a^*)}{b^*-a^*},
\end{equation}
as required.
\end{proof}

\begin{example}
Suppose that $X=\RR$ and that $C=[a,b]$, where
$0<a<b<\pinf$.
In each of the following items,
suppose that $f$ is as in the corresponding
item of Example~\ref{:ex:Legendre}.
Denote the corresponding
right and left Chebyshev centers by $x$ and $y$, respectively.
Then the following hold.
\begin{enumerate}
\item $\displaystyle x = y = \frac{a+b}{2}$. \\[+1mm]
\item $\displaystyle x = \frac{b-a}{\ln(b)-\ln(a)}$ ~and~
$\displaystyle y = \exp\Big(\frac{b\ln(b)-b-a\ln(a)+a}{b-a}\Big)$.\\[+1mm]
\item $\displaystyle x = \frac{ab\big(\ln(b)-\ln(a)\big)}{b-a}$ ~and~
$\displaystyle y = \frac{b-a}{\ln(b)-\ln(a)}$.
\end{enumerate}
\end{example}
\begin{proof}
This follows from Theorem~\ref{:t:kuehl}.
\end{proof}

\section{Generalizations and Variants}
\label{:sec:Shawn}

Chebyshev set and Klee set problems can be generalized to problems
involving functions.

Throughout this section,
\boxedeqn{
g\colon X\to \RRX \text{~~is lower semicontinuous and proper.}
} 
For convenience, we also set
\boxedeqn{
q := \tfrac{1}{2}\|\cdot\|^2.
}
\noindent
Recall that the \emph{Moreau envelope}
$e_{\lambda}g\colon X\to \RXX$
and the set-valued \emph{proximal mapping}
$P_{\lambda}g\colon X\rightrightarrows X$
are given by
\begin{equation}
x\mapsto e_{\lambda}g(x):=\inf_{w}\Big(g(w)+\frac{1}{2\lambda}\|x-w\|^2
\Big)
\end{equation}
and
\begin{equation}
x\mapsto P_{\lambda}g(x):=\argmin_{w}
\Big(g(w)+\frac{1}{2\lambda}\|x-w\|^2\Big).
\end{equation}
It is natural to ask: If $P_{\lambda}g$ is single-valued everywhere
on $\RR^n$, what can we say about the function $g$?

Similarly, define
$\phi_{\mu}g\colon X\to\RRX$ and
$Q_{\mu}g\colon X\rightrightarrows X$ by
\begin{equation}
y\mapsto \phi_{\mu}g(y):=\sup_{x}\Big(\frac{1}{2\mu}\|y-x\|^2-g(x)\Big),
\end{equation}
and
\begin{equation}
y\mapsto Q_{\mu}g(y):=\argmax_{x}\Big(\frac{1}{2\mu}\|y-x\|^2-g(x)
\Big).
\end{equation}
Again, it is
natural to ask: If $Q_{\mu}g$ is single-valued everywhere on $X$,
what can we say about the function $g$?
When $g=\iota_{C}$, then $P_{\lambda}g=P_{C}, Q_{\mu}g=Q_{C}$,
and we recover the classical Chebyshev and Klee set problems.

\begin{definition}
\
\begin{enumerate}
\item
The function $g$ is \emph{prox-bounded}
if there exists $\lambda>0$ such that
$e_{\lambda}g\not\equiv -\infty$.
The supremum of the set of all such $\lambda$
 is the threshold $\lambda_{g}$ of the prox-boundedness for $g$.
\item
The constant $\mu_{g}$ is defined to be
the infimum of all $\mu>0$ such that
$g-\mu^{-1}q$ is bounded below on $X$; equivalently,
$\phi_{\mu}g(0)<+\infty$.
\end{enumerate}
\end{definition}

\begin{fact}
\emph{(See \cite[Example 5.23, Example 10.32]{Rock98}.)}
Suppose that $g$ is prox-bounded with threshold $\lambda_{g}$, and
let $\lambda \in \left]0,\lambda_{g}\right[$.
Then $P_{\lambda}g$ is everywhere upper semicontinuous and locally
bounded on $X$, and $e_{\lambda}g$ is locally Lipschitz on $X$.
\end{fact}

\begin{fact}
\label{:theboss}
\emph{(See \cite[Proposition 4.3]{wang}.)}
Suppose that $\mu>\mu_{g}$.
Then $Q_{\mu}g$ is upper semicontinuous
and locally bounded on $X$, and
$\phi_{\mu}g$ is locally Lipschitz on $X$.
\end{fact}

\begin{definition}
\
\begin{enumerate}
\item We say that $g$ is \emph{$\lambda$-Chebyshev} if
$P_{\lambda}g$ is single-valued on $X$.
\item We say that $g$ is \emph{$\mu$-Klee} if $Q_{\mu}g$
is single-valued on $X$.
\end{enumerate}
\end{definition}

Facts~\ref{:hockey1} and \ref{:hockey3}
below concern Chebyshev functions and Klee functions;
see \cite{wang} for proofs.

\begin{fact}[single-valued proximal mappings]
\label{:hockey1}
Suppose that $g$ is prox-bounded with threshold $\lambda_{g}$, and let
$\lambda\in \left]0,\lambda_{g}\right[$.
Then the following are equivalent.
\begin{enumerate} \item
$e_{\lambda}g$ is continuously differentiable on $X$.\\[-1mm]
\item
$g$ is $\lambda$-Chebyshev, i.e.,
$P_{\lambda}g$ is single-valued everywhere.\\[-1mm]
\item $g+\lambda^{-1} q$ is essentially strictly convex.
\end{enumerate}
If any of these conditions holds, then
\begin{equation}
\label{:qiduan}
\nabla \big((g+\lambda^{-1}q)^{*}\big)=P_{\lambda}g \circ (\lambda \Id).
\end{equation}
\end{fact}

\begin{corollary}
\label{:hockey2}
The function $g$ is convex if and only if
$\lambda_{g}=+\infty$ and $P_{\lambda}g$ is single-valued on $X$ for every
$\lambda>0$.
\end{corollary}

\begin{fact}
[single-valued farthest mappings]\label{:hockey3}
Suppose that $\mu>\mu_{g}$. Then the following are equivalent.
\begin{enumerate}
\item $\phi_{\mu}g$ is (continuously) differentiable  on $X$.\\[-1mm]
\item
$g$ is $\mu$-Klee, i.e.,
$Q_{\mu}g$ is single-valued everywhere.\\[-1mm]
\item $g-\mu^{-1}q$ is essentially strictly convex.
\end{enumerate}
If any of these conditions holds, then
\begin{equation}
\label{:qiduan1}
\nabla \big((g-\mu^{-1}q)^*\big)=Q_{\mu}g(-\mu\Id).
\end{equation}
\end{fact}

\begin{corollary}\label{:hockey4}
Suppose that $g$ has bounded domain. Then
$\dom g$ is a singleton if and only if
for all $\mu>0$, the farthest operator $Q_{\mu}g$ is single-valued
on $X$.
\end{corollary}

\begin{definition}[Chebyshev points]
The set of $\mu$-Chebyshev points of $g$ is
$\argmin \phi_{\mu}g$.
If $\argmin \phi_{\mu}g$ is a singleton,
then we denote its unique element
by $p_{\mu}$ and we refer to $p_\mu$ as the \emph{$\mu$-Chebyshev point of
$g$}.
\end{definition}

The following result is new.

\begin{theorem}[Chebyshev point of a function]
Suppose that $\mu >\mu_{g}$.
Then the set of $\mu$-Chebyshev points is a singleton, and
the $\mu$-Chebyshev point is characterized by
\begin{equation}
p_\mu \in \conv Q_{\mu}g(p_\mu).
\end{equation}
\end{theorem}

\begin{proof}
As $\mu>\mu_{g}$, Fact~\ref{:theboss} implies that
\begin{equation}
y\mapsto \phi_{\mu}g(y)=
\frac{1}{2\mu}\|y\|^2+\left(-\frac{1}{\mu}q+g\right)^*(-y/\mu),
\end{equation}
is finite.
Hence $\phi_{\mu}g$ is strictly convex and super-coercive; thus,
$\phi_{\mu}g$ has a unique minimizer.
Furthermore, we have
\begin{equation}
\partial\phi_{\mu}g (y)=\frac{1}{\mu}
\big(y-\conv Q_{\mu}g(y)\big)
\end{equation}
by the Ioffe-Tikhomirov Theorem \cite[Theorem 2.4.18]{Zali02}.
Therefore,
\begin{equation}
0\in \partial\phi_{\mu}g (y) \quad \Leftrightarrow
\quad y\in \conv Q_{\mu}g(y),
\end{equation}
which yields the result.
\end{proof}

We now provide three examples to illustrate the
Chebyshev point of functions.

\begin{example}
Suppose that $g=q$. Then $\mu_{g}=1$ and for $\mu>1$, we have
\begin{equation}
\phi_{\mu}g\colon y\mapsto
\sup_{x}\bigg(\frac{1}{2\mu}(y-x)^2-\frac{x^2}{2}\bigg)=
\frac{y^2}{2(\mu-1)}.
\end{equation}
Hence the $\mu$-Chebyshev point of $g$ is
$p_{\mu}=0$.
\end{example}

\begin{example}
Suppose that $g=\iota_{[a,b]}$, where $a<b$.
Then $\mu_{g}=0$ and for  $\mu>0$, we have
\begin{equation}
\phi_{\mu}g \colon y
\mapsto \sup_{x}\bigg(\frac{1}{2\mu}(y-x)^2-\iota_{[a,b]}(x)\bigg)
=\begin{cases}
\frac{(y-b)^2}{2\mu} & \text{ if $y\leq \frac{a+b}{2}$,}\\
\frac{(y-a)^2}{2\mu} & \text{ if $y>\frac{a+b}{2}$.}
\end{cases}
\end{equation}
Hence $p_{\mu}=\frac{a+b}{2}$.
\end{example}

\begin{example} Let $a<b$ and
suppose that $g$ is given by
\begin{equation}
x\mapsto \begin{cases}
0 & \text{ if $a\leq x\leq \frac{a+b}{2}$},\\
1 & \text{ if $\frac{a+b}{2}< x\leq b$},\\
+\infty & \text{ otherwise}.
\end{cases}
\end{equation}
Then $\mu_{g}=0$, and when $\mu>0$ we have
\begin{align*}
\phi_{\mu}g(y) &= \sup_{x}\bigg(\frac{1}{2\mu}(y-x)^2 -g(x)\bigg)\\
&=\sup_{x}\begin{cases}
\frac{1}{2\mu}(y-x)^2 &\text{ if $a\leq x\leq \frac{a+b}{2}$}\\
\frac{1}{2\mu}(y-x)^2-1 &\text{ if $\frac{a+b}{2}<x\leq b$}\\
-\infty &\text{ otherwise}\\
\end{cases}
\\
&=\max\left\{\frac{(y-a)^2}{2\mu}, \frac{(y-(a+b)/2)^2}{2\mu},\frac{(y-b)^2}{2\mu}-1\right\},
\end{align*}
by using the fact that a strictly convex function
only achieves its maximum at the extreme points of its domain.
Elementary yet tedious calculations yield the following.
When $\mu>(a-b)^2/4$, we have
$$\phi_{\mu}g(y)=\begin{cases}
\frac{(y-b)^2}{2\mu}-1 & \text{ if $y<\frac{2\mu}{a-b}+\frac{a+3b}{4}$}\\
\frac{(y-(a+b)/2)^2}{2\mu} &\text{ if $\frac{2\mu}{a-b}+\frac{a+3b}{4}\leq y<\frac{3a+b}{4}$}\\
\frac{(y-a)^2}{2\mu} &\text{ if $y>\frac{3a+b}{4}$};
\end{cases}
$$
while when $0<\mu\leq (a-b)^2/4$, one obtains
$$\phi_{\mu}g(y)=
\begin{cases}
\frac{(y-b)^2}{2\mu}-1 & \text{ if $y<\frac{\mu}{a-b}+\frac{a+b}{2}$}\\
\frac{(y-a)^2}{2\mu} & \text{ if $y\geq \frac{\mu}{a-b}+\frac{a+b}{2}$}.
\end{cases}
$$
Hence, the Chebyshev point of $g$ is
\begin{equation*}
p_\mu = \begin{cases}
\displaystyle \frac{3a+b}{4}, & \text{if $\mu> (a-b)^2/4$;}\\[+4mm]
\displaystyle \frac{\mu}{a-b}+\frac{a+b}{2}, &\text{if $0<\mu \leq (a-b)^2/4$.}
\end{cases}
\end{equation*}
\end{example}

\section{List of Open Problems}

\label{:sec:op}

\begin{description}
\item[\textbf{Problem 1.}]
Is the assumption that $f$ be supercoercive in
Fact~\ref{:f:leftCheb} really essential?
\item[\textbf{Problem 2.}]
Are the assumptions that $f$ have full domain and
that $\overline{C^*}\subseteq U^*$ in
Fact~\ref{:f:rightCheb} really essential?
\item[\textbf{Problem 3.}]
Does there exist a Legendre function $f$ with full domain
such that $f$ is not quadratic yet every nonempty closed convex subset of
$X$ is $\fD{}$-Chebyshev?
In view of Fact~\ref{:f:leftCheb},
the gradient operator $\nabla f$ of such a function would be nonaffine
and it would preserve convexity.
\item[\textbf{Problem 4.}]
Is it possible to characterize the class of
$\fD{}$-Chebyshev subsets of the
strictly positive orthant when $f$ is the negative entropy?
Fact~\ref{:f:Noll} and Example~\ref{:ex:curious} imply
that this class contains not only all closed convex
but also some nonconvex subsets.
\item[\textbf{Problem 5.}]
Is the assumption that $f$ be supercoercive in Fact~\ref{:f:smallbP}
really essential?
\item[\textbf{Problem 6.}]
Is the assumption that $f$ be supercoercive in
Fact~\ref{:f:smallbQ} really essential?
{
\item[\textbf{Problem 7.}]
For the Chebyshev functions and Klee functions, we have used the halved Euclidean distance. What are characterizations
of $f$ and Chebyshev point of $f$ when one uses the Bregman distances?

\item[\textbf{Problem 8.}]
How do the results on Chebyshev functions and Klee functions
extend to Hilbert spaces or even general Banach spaces?
}
\end{description}

\section{Conclusion}

Chebyshev sets, Klee sets, and Chebyshev centers are well known
notions in classical Euclidean geometry. These notions have been studied
traditionally also in infinite-dimensional setting or with respect to
metric distances induced by different norms. Recently, a new framework was
provided by measuring the discrepancy between points differently,
namely by Bregman distances, and
new results have been obtained that generalize the classical
results formulated in Euclidean spaces. These results are fairly well
understood for Klee sets and Chebyshev centers with respect to Bregman
distances; however, the situation is much less clear for Chebyshev sets.

The current state of the art is reviewed in this paper and several new
results have been presented. The authors
hope that in the list of open problems (in Section~\ref{:sec:op})
will entice the reader to make further progress on this fascinating topic.

\section*{Acknowledgments}
Heinz Bauschke was partially supported by the Natural Sciences and
Engineering Research Council of Canada and by the Canada Research Chair
Program.
Xianfu Wang was partially supported by the Natural Sciences and Engineering
Research Council of Canada.
%

\end{document}